\documentclass[11pt]{article}

\usepackage{amssymb,amsmath}
\usepackage{color}
\newtheorem{Theorem}{Theorem}[section]

\newtheorem{Lemma}{Lemma}[section]
\newtheorem{Definition}{Definition}[section]
\newtheorem{Corollary}{Corollary}[section]

\numberwithin{equation}{section}

\newcommand{\non}{\nonumber}
\newcommand{\cR}{\mathbb R}
\newcommand{\cN}{{\mathbb N}}

\newcommand{\eq}[1]{\mbox{\rm {(\ref{#1})}}}

\begin{document}
\title{\Large\bf Weak solutions to an initial-boundary value problem
for a continuum equation of  motion of grain boundaries
}
\author{{\small\sc  Peicheng Zhu$^1$\footnote{
E-mail: pczhu@shu.edu.cn}, {\small\sc Lei Yu$^1$\footnote{E-mail: yulei@shu.edu.cn} }
and
Yang Xiang$^2$\footnote{E-mail: maxiang@ust.hk} }\\
{\small\sc $^1$Department of   Mathematics, Shanghai University,}\\
{\small\sc  Shanghai 200444, P.R. China} \\
{\small\sc $^2$Department of Mathematics, Hong Kong University of  Science and Technology,}\\
{\small\sc Clear water bay, Kowloon, Hong Kong  }
}

\date{}

\maketitle

\noindent{\bf Abstract.} We investigate
 an initial-(periodic-)boundary value problem for a continuum equation,
which is a model for motion of grain boundaries based
on the underlying microscopic mechanisms of line defects (disconnections) and integrated the effects of a
diverse range of thermodynamic driving forces. We first prove the
 global-in-time existence and uniqueness of weak solution to this initial-boundary value problem in
 the case with positive equilibrium disconnection density parameter $B$, and then investigate the
 asymptotic behavior of  the solutions  as $B$ goes to zero.
The main difficulties in the proof of main theorems are due to the degeneracy of $B=0$,
a non-local term with singularity,
and a non-smooth coefficient of the highest derivative associated with the gradient of the unknown.
The key ingredients in the proof are the energy method, an estimate for
a singular integral of the Hilbert type,  and a compactness lemma.

\medskip
\noindent{\bf Keywords.} Motion of grain boundaries;
Initial-boundary value problem; Global existence; Weak solutions; Disconnections

%
%

\section{Introduction}

 A polycrystalline material can be regarded as a network
 of grain boundaries (GBs) on the mesoscale.
This GB network has a great impact on a wide range of
 materials properties, such as strength, toughness,
 electrical conductivity, and its evolution is
 important for engineering materials~\cite{SB95}. Grain boundaries are
 the interfaces between  differently oriented crystalline grains, which
 are a kind of two-dimensional defects in
 materials. Grain boundary migration controls many
 microstructural evolution processes in materials. Since GBs are interfaces between crystals, the microscopic
mechanisms by which they move are intrinsically different
from other classes of interfaces, such as solid-liquid interfaces and
 biological cell membranes.

Recent experiments and atomistic simulations have shown that the microscopic mechanism of GB migration is associated
with the motion of topological line defects, i.e., disconnections~\cite{A72,KS80,HP96,MTP04,HPL07,RMLC13,RLCMM13,TS17}.
This  dependence on microscopic structures enables broad-range and deep understandings of GB migration, e.g, the stress-driven motion and the shear coupling effect~\cite{LEWP53,CT04}, which cannot be described by the classical motion by mean curvature models (driven by capillary forces)~\cite{SB95}.

A new continuum equation  for motion of grain boundaries based on the  underlying disconnection mechanisms was developed by Zhang {\it et al.} \cite{13} in 2017. This continuum model integrates the effects of a
diverse range of thermodynamic driving forces including the stress-driven motion and is able to describe the shear coupling effect during the GB motion.
Generalizations of this continuum model with multiple disconnection modes and GB triple junctions have been futher developed~\cite{WX19,WX20,ZX21}.



In the present article, we will study the existence of weak solutions to the
initial-boundary value problem of the continuum equation for GB motion developed in Ref.~\cite{13}, which reads
\begin{eqnarray}
 h_{t} = - M_{d}\big((\sigma_{i}+\tau)b+\Psi H-\gamma H h_{xx}\big)(|h_{x}|+B)
 \label{1.1}
\end{eqnarray}
for $(t,x)\in(0,\infty)\times\Omega$,  where $\Omega=(a,d)$.
The boundary and initial conditions are
\begin{eqnarray}
 h |_{x=a} &=& h |_{x=d},\quad h_{x} |_{x=a} = h_{x} |_{x=d},\ (t,x)\in(0,T_{e})\times\partial\Omega,
 \label{1.2}\\
 h(0,x) &=& h_{0}(x),\ x\in\Omega,
 \label{1.3}
\end{eqnarray}
where
\begin{eqnarray}
 \sigma_{i}(t,x) = {\rm P.V.}\int_{-\infty}^{\infty}\frac{K\beta h_{x}(t,x_{1})}{x-x_{1}}dx_{1},
 \label{1.4}
\end{eqnarray}
and
$$
 \beta=\frac{b}{H},\ K = \frac{\mu}{2\pi(1-\nu)},\
  B = \frac{2H}{a} e^{-F_{d}/(k_{B}T)}.
$$

The unknown function $h$ in Equation \eqref{1.1} is
 the height of grain boundary from reference line,
and $\sigma_{i}(x,t)$ is the stress due to the elastic interaction between disconnections based on their dislocation nature~\cite{HL,Xiang2006}.
The parameters
$b$ and $H$ respectively are the Burgers vector and step height of a disconnection, $\mu$ and $\nu$ respectively are the shear modulus and Poisson ratio,
$\gamma$ is the GB energy, $\tau$ is the applied stress, $\Psi$ is energy jump across the GB, and $M_{d}$ is the mobility constant. The parameter $B$ is associated with the equilibrium density of the disconnection pairs, where $F_{d}$ is the disconnection formation energy, $a$ is the lattice constant,  $k_{B}$ is the Boltzmann constant, $T$ is the temperature, and $\frac{1}{a} e^{-F_{d}/(k_{B}T)}$ is the equilibrium disconnection density.

 We will study the existence of weak solutions for both cases of $B>0$ and $B=0$.
Note that the regime of $B\rightarrow 0$ means small equilibrium disconnection density or large slope of the grain boundary profile, and  when $B=0$, the equation \eqref{1.1} is degenerate at those points where $h_x=0$. Numerical results in Ref.~\cite{13} showed that sharp corners may be developed in the GB profile in the case of $B=0$.


The  difficulties in the proofs of the existence and uniqueness theorems come from
the non-local term with singularity
together with a non-smooth coefficient associated with $|h_x|$ of the highest derivative $h_{xx}$, and the degeneracy of the equation in the case of $B=0$. To estimate this singular integral term, we employ a theorem
in the book by Stein~\cite{18}. Regularization is performed so that the coefficient of the $h_{xx}$ term is smooth and uniformly bounded from below, and then  compactness lemmas are employed to obtain the results for the original equations. Note that dependence  on non-smooth gradient terms in the coefficient of the highest derivative also appeared in the phase field models proposed by Alber and Zhu  in~\cite{1,3} to describe the evolution of an interface driven by configurational forces, and properties of solutions have been obtained~\cite{2,4,5,15,16z,16}.
Models with  non-smooth gradient terms have also been investigated by Acharya, et al~\cite{AMZ10}
and Hilderbrand, et al~\cite{H10}.

\subsection{Interpretation of the formula for $\sigma_i$}

This subsection is intended to give an explanation of the formula for $\sigma_i$.
The material we consider is
 normally in a bounded domain and we assume the spatial periodic boundary conditions,
 which means that the unknown
 $h$ is defined over $\Omega$; however in  formula \eqref{1.4},   the integral   domain
  is   $\cR=(-\infty,+\infty)$, which implies $h$ should be defined over
  $\cR$.
  Letting $L = d-a$ be the smallest positive period, choosing $x\in\Omega$ we then arrive at
\begin{eqnarray}
 \sigma_{i}(t,x)
 &=&  {\rm P.V.}\int_{-\infty}^{\infty}\frac{K\beta h_{x}(t,x_{1})}{x-x_{1}}dx_{1}
  \nonumber\\
 &=&  \sum_{k\in \mathbb{Z}}  {\rm P.V.}\int_{a+kL}^{a+(k+1)L}\frac{K\beta h_{x}(t,x_{1})}{x-x_{1}}dx_{1}
  \nonumber\\
 &=& \sum_{k\in \mathbb{Z}}  {\rm P.V.}\int_{a}^{d}\frac{K\beta h_{x}(t,y)}{x-y+kL}dy\nonumber\\
 &=:&  \sigma_{i1} + \sigma_{i2} ,
 \label{1.6}
\end{eqnarray}
where
\begin{eqnarray}
 \sigma_{i1} &=:&  {\rm P.V.}\sum_{k\in \mathbb{Z}\backslash\{0\}}\int_{a}^{d}\frac{K\beta h_{x}(t,y)}{x-y+kL}dy
 ,
  \label{1.6a}\\
 \sigma_{i2} &=:&  {\rm P.V.} \int_{a}^{d}\frac{K\beta h_{x}(t,y)}{x-y }dy .
 \label{1.6b}
\end{eqnarray}
Here $\mathbb{Z}$ denotes the set of integers. Observing that for  $k\not=0$
 we see that the function $\frac{1}{x-y+kL}$ is positive and monotonically increasing
  in $y$ for any fixed $x$ and $L$, thus by applying the second mean value theorem of
   integrals, we conclude that  there exists a number $\eta\in [a,d]$ such that
$$
 \int_{a}^{d}\frac{K\beta h_{x}(t,y)}{x-y+kL}dy = \frac{ 1}{x-d+kL}\int_{\eta}^{d} {K\beta h_{x}(t,y)} dy ,
$$
while for the case that  $k = 0$,   $\int_{a}^{d}\frac{K\beta h_{x}(t,y)}{x-y }dy $ is a singular integral for general $h$,
from which one thus finds that for many kinds of $h$, the
series in $\sigma_i$ may diverge.

Therefore one must understand  Principal Value in the formula of $\sigma_i$ both for
a series $\sigma_{i1}$ and for singular integral $\sigma_{i2}$, more precisely, the term
\begin{eqnarray}
 \sigma_{i1} = \sum_{k=1}^\infty
 \int_{a}^{d}\left(\frac{K\beta h_{x}(t,y)}{x-y+kL} + \frac{K\beta h_{x}(t,y)}{x-y-kL}\right) dy ,
 \label{1.6c}
\end{eqnarray}
which implies that
\begin{eqnarray}
 |\sigma_{i1}| \le C \|h_{x}\|_{L^2(\Omega)} \sum_{k=1}^\infty \frac1{(kL)^2} \le C \|h_{x}\|_{L^2(\Omega)},
 \label{1.6d}
\end{eqnarray}
and
$$
 \sigma_{i2} =\lim_{\varepsilon\to 0} \int_{\{|x-y|> \varepsilon\}\cap \Omega}
 \frac{K\beta h_{x}(t,y)}{x-y }dy,
$$
whose bound can be evaluated by employing Theorem~\ref{thm2.2} suppose that $h_x\in L^2(\Omega)$.

{\bf Remark.} We would like to point out  another  way to interpret the singular integral $\sigma_i$ with the following modification:
\begin{eqnarray}
 \sigma_{i}(t,x)
 &=&  {\rm P.V.}\int_{-\infty}^{\infty}\frac{K\beta h_{x}(t,x_{1})}{ D(x-x_{1}) }dx_{1} =:  \sigma_{i1} + \sigma_{i2},
 \label{1.7a}
\end{eqnarray}
where $D(x-x_{1})$ is defined by
\begin{eqnarray}
 D(x-x_{1}) = \left\{
 \begin{aligned}
  & {\rm sgn}(x_{1}-x)|x-x_{1}|^{1+\zeta}, & {\rm if}\ x\not\in [a,d];\\
  & x-x_{1}, & {\rm if}\  x\in [a,d]
 \end{aligned}
  \right.
 \label{zeta}
\end{eqnarray}
with an arbitrarily given positive constant $\zeta$.
We can then conclude that the series $\sigma_{i1}$ converges and the singular integral $\sigma_{i2}$ may
be treated as in the previous method. We will not use this method in the proofs in this paper.
\subsection{Main results}

We first perform nondimensionalization of the quation. Using $M_d\mu$ as the time unit, $\mu$ the unit of $\sigma_i$, $\tau$ and $\Psi$, $L_0$ the unit of the length scale of the continuum equation, and $\mu L_0$ the unit of $\gamma$, we have the dimensionless form of the equation.
Further introducing parameters
$$
 \alpha_{1} = \gamma H,\
 \alpha_{2} =  b, \
 \alpha_{3} =  \tau b+\Psi H,
$$
where all the quantities are in dimensionless form, equation  \eqref{1.1} can be written as
\begin{eqnarray}
 h_{t} - \frac{\alpha_{1}}{2}\left( |h_{x}|h_{x}  + 2Bh_{x}\right)_x
 + (\alpha_{2}\sigma_i + \alpha_{3} )(|h_{x}| + B) =0.
 \label{1.9}
\end{eqnarray}
Here we have used the  formula
 $ (|y|y)^{'} = 2|y|$. From now on, we will use this nondimensionalized equation with the dimensionless parameters described above.

\vskip0.25cm
To define weak solutions to the initial-boundary value
 problem  (\ref{1.1}) -- (\ref{1.3}), we denote by $\Omega=(a,d)$
    a bounded open interval with constants
  $a<d$,  by $T_{e}>0$   an arbitrary constant, and by
   $Q_{T_e}$ the domain $(0,T_{e})\times\Omega$.  Define
$$
 (v_1,v_2)_{Z}=\int _{Z}v_1(y)v_2(y)dy
$$
for $Z=\Omega$ or $Z=Q_{T_{e}}$. Moreover, if $v$ is a function
defined on $Q_{T_{e}}$, we use $v(t)$ to represent the mapping
$x\mapsto v(t,x)$ and sometimes write $ v = v(t)$ for convenience.

\vskip0.2cm
\noindent{\bf Statement of the main results.} Our main results are concerned
with the existence and uniqueness of weak solution to an initial-boundary
 value problem.

\begin{Definition}\label{def1.1}
 Let $h_{0}\in L^{1}(\Omega)$. A function $h$ with
\begin{eqnarray}
 h\in  L^{2}(0,T_{e};H^{1}_{{\rm per}}(\Omega))
 \label{1.10}
\end{eqnarray}
is called a weak solution to  problem \eqref{1.1} -- \eqref{1.3}, if for
all $\varphi\in C_{0}^{\infty}((-\infty,T_{e})\times\Omega)$, there holds
\begin{equation}
 (h,\varphi_{t})_{Q_{T_{e}}} - \frac{\alpha_{1}}{2} ( |h_{x}|h_{x} + 2B h_{x},\varphi_{x} )_{Q_{T_{e}}}
 - ((\alpha_{2}\sigma_i + \alpha_{3} )(|h_{x}| + B)  , \varphi) _{Q_{T_{e}}}
 + (h_{0},\varphi(0))_{\Omega}=0.
 \label{1.11}
\end{equation}

\end{Definition}

We then have
\begin{Theorem}\label{thm1.1} Assume that $\gamma H $ is sufficiently greater than $b$,
 and $h_{0}\in H^{1}_{{\rm per}}(\Omega)$. Then there exists a unique
  weak solution $h$ to
 problem \eqref{1.1} -- \eqref{1.3} with $B>0$, which in addition
to \eqref{1.10}, satisfies
\begin{eqnarray}
 h\in L^{\infty}(0,T_{e};H^{1}_{{\rm per}}(\Omega)),\ \ h_{x}\in L^{2}(0,T_{e};H^{1}_{{\rm per}}(\Omega))
 \cap L^{{3}}(Q_{T_{e}}),
 \label{1.12a} \\
 h_{t}\in L^{\frac{4}{3}}(Q_{T_{e}}),\ \
 |h_{x}|h_{x} \in L^{\frac{4}{3}}(0,{T_{e}}; W^{1,\frac{4}{3}}_{{\rm per}}(\Omega)).
 \label{1.12}
\end{eqnarray}

\end{Theorem}

We are also interested in the limit as $B\to 0$.
\begin{Definition}\label{def5.1}
 Let $h_{0}\in L^{1}(\Omega)$. A function $h$ with
\begin{eqnarray}
 h\in  L^{2}(0,T_{e};H^{1}_{{\rm per}}(\Omega))
 \label{1.10B}
\end{eqnarray}
is called a weak solution to  problem \eqref{1.1} -- \eqref{1.3} with $B=0$, if for
all $\varphi\in C_{0}^{\infty}((-\infty,T_{e})\times\Omega)$, there holds
\begin{eqnarray}
 (h,\varphi_{t})_{Q_{T_{e}}} - \frac{\alpha_{1}}{2} ( |h_{x}|h_{x} ,\varphi_{x} )_{Q_{T_{e}}}
 - ((\alpha_{2}\sigma_i + \alpha_{3} ) |h_{x}|, \varphi) _{Q_{T_{e}}}
 +  (h_{0},\varphi(0))_{\Omega}=0.
 \label{1.11B}
\end{eqnarray}

\end{Definition}

We denote a solution to  problem \eqref{1.1} -- \eqref{1.3} by $h_B$, then $h_B$
converges $h$ almost everywhere $(t,x)$ over $Q_{T_e}$, and $h$ satisfies \eq{1.11B}.
\begin{Theorem}\label{thm1.1B} Assume that $\gamma H $ is sufficiently greater than $b$,
 and $h_{0}\in H^{1}_{{\rm per}}(\Omega)$. Then there exists a weak solution $h$ to
 problem \eqref{1.1} -- \eqref{1.3} with $B=0$, which in addition
to \eqref{1.10B}, satisfies
\begin{eqnarray}
 h\in L^{\infty}(0,T_{e};H^{1}_{{\rm per}}(\Omega)),\ \ h_{x}\in  L^{{3}}(Q_{T_{e}}),
 \label{1.12aB} \\
 h_{t}\in L^{\frac{4}{3}}(Q_{T_{e}}),\ \
 |h_{x}|h_{x} \in L^{\frac{4}{3}}(0,{T_{e}}; W^{1,\frac{4}{3}}_{{\rm per}}(\Omega)),\\
 (|h_{x}|h_{x})_t \in L^{1}(0,{T_{e}}; H^{-2}_{{\rm per}}(\Omega)).
 \label{1.12B}
\end{eqnarray}

\end{Theorem}

\noindent{\bf Remarks.} {\it 1. In the original units,  the assumption
that $\gamma H $ is sufficiently greater than $b$ means that $\gamma\gg (b/H)\mu L_0$, where $L_0$ is the length scale of the continuum equation.

2. For the regularity of the solution $h$, we have a more regular weak solution in the case of $B>0$
than  that in the case of $B=0$.  This result agrees with the numerical results obtained in Ref.~\cite{13} in which sharp corners were developed in $h$ in the case of $B=0$.
}

\vskip0.2cm
\noindent{\bf Notations.}
$C, C(\cdot)$ denote, respectively,  universal constants which may vary from line to line.
 and $C(\cdot)$ depends on its argument(s).
Greek letters $\varepsilon,\ \zeta$ are small positive numbers which are normally
 assumed to be small.
$\kappa$ is taken in $(0,1]$, which will be sent to zero.
${T_e}$ (or $t_e$) denotes a positive constant related to time, the life of a solution.

 Let $p,\ q$ be  real numbers such that $p,\ q\ge 1$.
 Let $\cN$ be the set of natural number and $\cN_+ = \cN\cup \{0\}$, and $\cR^d$ be $d$-dimensional
 Euclidean space.

 $\Omega$ denotes an open, bounded, simple-connected domain in $\cR^d$ with natural number
 $d$, with smooth boundary $\partial \Omega$. It represents
the material points of a solid body. $Q_t=(0,t)\times
\Omega$, and its parabolic boundary ${\cal P}Q_t$ is defined by $ {\cal P}Q_t:=(\partial\Omega\times [0,t))\cup(\Omega\times\{0\})$.

 $L^p(\Omega)$ are the Sobolev spaces of $p$-integrable
real functions over $\Omega$ endowed with the norm
$$
 \|f\|_{L^p(\Omega)} = \left(\int_\Omega|f(x)|^pdx\right)^\frac1p, \mbox{ if }p<\infty;
 \ \|f\|_{L^\infty(\Omega)} = \mbox{ ess} \sup_{x\in\Omega}|f(x)|.
$$
 Throughout this article, the norm of $L^2(\Omega)$ is denoted by $\|\cdot\|$,
 and the norm of $L^2(Q_{T_{e}} )$ is denoted by $\| \cdot \|_{Q_{T_{e}}}$.

Let $\Omega$ be an $n$-dimensional cuboid. Let $\alpha\in \cN_0^d$  be a multi-index and $|\alpha|$ be its length, where
 $\cN_0=\cN\cup\{0\}$.
$D^\alpha f$ is the $|\alpha|$-th order weak derivatives. Define the space
 $W^{m,p}_{{\rm per}}(\Omega)=\{ f\in L^p(\Omega)\mid D^\alpha f\in L^p(\Omega) \mbox{ for all } \alpha
\mbox{ such that }
|\alpha|\le m,  \mbox{ and } \gamma_jf|_{{\rm on\ one\ face}} = (-1)^j \gamma_jf|_{{\rm on\ the\ corresponding\ face}},
\mbox{ for } j=0,1,\cdots,m-1\}$ endowed with norm
$$
 \|f\|_{W^{m,p}_{{\rm per}}} = \left( \sum_{|\alpha|\le m} \|D^\alpha f\|_{L^p(\Omega)}^p \right)^\frac1p,
$$
where $\gamma_j$ are the trace operators. And
$W_0^{m,p}(\Omega)$ is the closure of $C_0^\infty(\Omega)$ in   the norm $\|\cdot\|_{W^{m,p}(\Omega)}$.
For $p=2$, $H^m_{{\rm per}} (\Omega):= W^{m,2}_{{\rm per}}(\Omega)$,
 $H^{-m}_{{\rm per}} (\Omega)$ denotes the dual space of $H^m_{{\rm per}} (\Omega)$.
$H^m _0(\Omega):= W^{m,2}_0(\Omega)$.

Let $q,p\in \cR$ such that $q,p\ge 1$.
\begin{eqnarray}
 L^q(0,T;L^{p}(\Omega))
  &:=&\Big\{f \mid f \mbox{ is Lebesgue measurable such that }  \non\\
  && \|f\|_{L^q(0,t;L^p(\Omega))} := \left(\int_0^t\left(\int_\Omega |f|^pdx\right)^{\frac{q}p}d\tau \right)^\frac1q<\infty \Big\},
 \non
\end{eqnarray}
and
$L^q(0,t;W^{m,p}_{{\rm per}}(\Omega)):=\{f\in L^q(0,t; L^{p}(\Omega))
\mid \int_0^t \|f(\cdot,\tau)\|_{W^{m,p}_{{\rm per}}} ^q d\tau<\infty \}$. See, e.g., \cite{14}.
 We also need some function spaces: For non-negative integers  ${m},  {n}$,
real number $\alpha\in(0,1)$ we denote by $C^{m+\alpha}(\overline{\Omega})$ the
space of \textit{m}-times differentiable functions on $\overline{\Omega}$, whose
$m$th derivative is H\"{o}lder continuous with exponent $\alpha$. The space
$C^{\alpha,\frac{\alpha}{2}}(\overline{Q}_{T_{e}})$ consists of all functions
on $\overline{Q}_{T_{e}}$, which are H\"{o}lder continuous in the parabolic distance
$$
 d((t,x),(s,y)):=\sqrt{|t-s|+|x-y|^{2}}
$$
$C^{m,n}(\overline{Q}_{T_{e}})$ and $C^{m+\alpha,n+\frac{\alpha}{2}}(\overline{Q}_{T_{e}})$,
respectively, are the spaces of functions, whose ${x}$-derivatives up to
order ${m}$ and ${t}$-derivatives up to order \textit{n} belong to
$C(\overline{Q}_{T_{e}})$ or to $C^{\alpha,\frac{\alpha}{2}}(\overline{Q}_{T_{e}})$,
respectively.

\vskip0.2cm
\noindent{\bf Organization of rest of this article.} The main results of this article are  Theorem~\ref{thm1.1}  and Theorem~\ref{thm1.1B}.
 The remaining sections  are devoted to the proofs of these theorems.
In Section~2, we construct an approximate initial-boundary value problem,
and prove by employing the Leray-Schauder fixed-point theorem,
the existence of classical solutions to
this problem. Then we derive in Section~3 {\it a prior} estimates which are uniform in a small
parameter, in which we used some results about singular integrals.
 In Section~4, by making use of the Aubin-Lions lemma and properties of
  strong convergence, we prove the existence and uniqueness of weak solution to the original initial-boundary  value problem when $B$ is a positive constant. Section~5 is
devoted to the study of asymptotic behavior as $B$ goes to zero, to this end we establish {\it a prior} estimates that are independent of $B$, consequently we
prove Theorem~\ref{thm1.1B}.

%

\section{Existence of solutions to the modified problem}

To prove Theorem~\ref{thm1.1}, we  construct the following approximate
 initial-boundary value  problem.
\begin{eqnarray}
   h_{t} -  {\alpha_{1}} \left( \int_{0}^{h_{x}}|p|_\kappa dp +  Bh_{x}\right)_x
 &+& (\alpha_{2}\sigma_i^\kappa + \alpha_{3} )(|h_{x}|_\kappa + B)
  = 0, \ {\rm in}\ Q_{T_e},
 \label{2.1}\\
 h |_{x=a} &=& h |_{x=d},\quad h_{x} |_{x=a} = h_{x} |_{x=d}, \ {\rm on}\  \partial\Omega\times [0,T_{e}],\ \
 \label{2.4}\\
 h(0,x)  &=& h^\kappa_{0}(x), \ {\rm in}\  \Omega.
 \label{2.5}
\end{eqnarray}
Here  $\kappa>0$ is a constant,  we used the notation
\begin{eqnarray}
 |p|_{\kappa} := \sqrt{|p|^{2}+\kappa^{2}}
 \label{2.2}
\end{eqnarray}%
to replace the function $|p|$ to smooth the coefficient of the principal term in
 \eq{2.1} and to guarantee  that equation is uniformly parabolic from below.
 And $\sigma_i^\kappa = \sigma_{i1}^\kappa  + \sigma_{i2}^\kappa $, where $\sigma_{i1}^\kappa$
 is obtained via replacing $h$ in $\sigma_{i1} $ by a function $h^\kappa\in L^2(0,T_e;H^2_{per}(\Omega))$, and
$\sigma_{i2}^\kappa $ is defined, for this $ h^\kappa$,  by
\begin{eqnarray}
 \sigma_{i2}^\kappa = K\beta  \int_{\{|y|>\kappa\}\cap \Omega }\frac{  h^\kappa_{x}(t,x-y)}{y}dy.
 \label{1.4a}
\end{eqnarray}
The initial data $h^\kappa_{0}(x)$ is chosen such that $h^\kappa_{0}\in C^\infty(\Omega)$ and
$$
 \|h^\kappa_{0} - h_{0}\|_{ H^1_{\rm per}(\Omega)} \to 0.
$$

We now state the existence of classical solution to problem \eq{2.1} -- \eq{2.5} as follows.
%
%
%
\begin{Theorem} \label{thm2.1}
Suppose that $\gamma H $ is sufficiently greater than $b$ and
 the initial data $h^\kappa_{0} $  satisfies the compatibility conditions
$ h_{0}^\kappa(a) = h_{0}^\kappa(d),\quad h^\kappa_{0x} (a) = h^\kappa_{0x} (d),
\quad h^\kappa_{0xx} (a) = h^\kappa_{0xx} (d)$, and
$ h^\kappa_{t} (0, a) + \alpha_{2}\sigma_i^\kappa(0,a) (|h^\kappa_{0x} (a)| + B)
  = h^\kappa_{t} (0, d) + \alpha_{2}\sigma_i^\kappa(0,d) (|h^\kappa_{0x} (d)| + B)$.

Then there exists a classical solution $h^\kappa$ to problem \eqref{2.1} -- \eqref{2.5}
 such that
$$
 h^\kappa_{xt} \in L^2(Q_{t_e}),\
 \|h^\kappa\|_{ C^{\beta/2,1+\beta}(\bar Q_{t_e}) } \le C_\kappa.
$$

\end{Theorem}

We now present the strategy of the proof of Theorem~\ref{thm2.1}.
Since there is a non-local term $\sigma_i^\kappa$ in Eq. \eqref{2.1}, we shall employ the Leray-Schauder
fixed-point theorem to prove this theorem. To this end, we first modify that equation as
\begin{eqnarray}
 h_{t} - \alpha_{1}h_{xx}(|h_{x}|_{\kappa}+B) &+&
   (\alpha_{2}\lambda \hat\sigma_i^\kappa + \alpha_{3} )(| h_{x}|_{\kappa}  + B)=0,
   \ {\rm in}\ Q_{T_e},
 \label{2.1a}\\
  h |_{x=a} &=& h |_{x=d},\quad  h_{x} |_{x=a} = h_{x} |_{x=d}, \ {\rm on}\  \partial\Omega\times [0,T_{e}],
 \label{2.4a}\\
 h(0,x) &=& \lambda h^\kappa_{0}(x), \ {\rm in}\  \Omega,
 \label{2.5a}
\end{eqnarray}
where $\lambda\in [0,1],\ \hat\sigma_i^\kappa = \hat\sigma_{i1}^\kappa + \hat\sigma_{i2}^\kappa$,
and $\hat\sigma_{i2}^\kappa$ is defined by
\begin{eqnarray}
 \hat\sigma_{i2}^\kappa = K\beta  \int_{\{|y|>\kappa \} \cap \Omega }\frac{\hat  h_{x}(t,x-y)}{y}dy.
 \label{1.4b}
\end{eqnarray}

We take $0 < \alpha < 1$ and define for any $\hat h \in {\cal B} := C^{\alpha/2 ,1+\alpha}(\bar Q_{t_e} )$
a mapping $P_\lambda : [0, 1] \times {\cal B} \to {\cal B}$;  $\hat h \mapsto  h$ where $h$ is the solution  to
problem \eq{2.1a} -- \eq{2.5a}, and the existence of solutions to this problem
can be found, e.g., in Theorem~4.1,~P. 558, Ladyshenskaya et al
 \cite{14}  with slight modifications.

Next we derive {\it a priori} estimates which may depend on the parameter $\kappa$. In the rest part  of this
section, we assume that the conditions in  Theorem~\ref{thm2.1} are met,
 and there exists a unique solution to \eq{2.1a} -- \eq{2.5a}, which means that   $\hat\sigma_i^\kappa$ in
\eq{2.1a} is replaced by $ \sigma_i^\kappa$.
\begin{Lemma} \label{lm-est-kappa-1} There holds for any $t\in [0,T_e]$ that
\begin{eqnarray}
 \|h^\kappa(t)\|^2_{H^1(\Omega)} + \int_0^t \|h^\kappa_{xx}(\tau)\|^2 d\tau \le C_\kappa.
 \label{est-kappa-1}
\end{eqnarray}
Here $C_\kappa$ is a constant which may depend on $\kappa$.

\end{Lemma}

 This estimate is easier to  obtain than those in Section~3 that are independent of $\kappa$,
 so we omit the details
  of the derivation for most of them. We also need the following  estimates.
\begin{Lemma} \label{lm-est-kappa-2} There holds for any $t\in [0,T_e]$ that
\begin{eqnarray}
 \|h^\kappa_{xx}(t)\|^2 + \|h^\kappa_t(t)\|^2  +  \int_0^t \|h^\kappa_{xt}(\tau)\|^2 d\tau \le C_\kappa.
 \label{est-kappa-2}
\end{eqnarray}

\end{Lemma}

This estimate is {\it not} necessary for the proof of the existence of weak solutions,
thus we  give the main idea on deriving it. For the sake of reader's convenience, we present
some tools as follows.  First we recall the Gronwall lemma.
\begin{Lemma}[Gronwall Lemma]\label{Gronwall}  Let $ y, A, B$ be functions satisfying
that $A(t), B(t)$  are integrable over $[0,t_e ]$ and $y(t) \ge  0$ is absolutely continuous
function. Then
$$
 y'(t) \le A(t)y(t) + B(t),\ for\ a.e.\ t,
$$
implies
$$
 y(t) \le  y(0){\rm e}^{\int_0^t A(\tau)d\tau} + \int_0^t B(s){\rm e}^{\int_s^t A(\tau)d\tau}ds.
$$

\end{Lemma}

\begin{Lemma}[Aubin-Lions] \label{thm2.3}
Let $B_{0}$ and $B_{2}$ be reflexive Banach spaces and let $B_{1}$
be a Banach space such that $B_{0}$ is compactly embedded in $B_{1}$
and that $B_{1}$ is embedded in $B_{2}$. For $1\leq p_{0}, p_{1}\le +\infty$,
define
$$
 W=\left\{f\mid f\in L^{p_{0}}(0,T_e;B_{0}),\  \frac{df}{dt}\in L^{p_{1}}(0,T_e;B_{2})\right\}.
$$
$($\textit{i}$)$ if $p_{0}<+\infty$, then the embedding of $W$ into
$L^{p_{0}}(0,T_e;B_{1})$ is compact.\\
$($\textit{ii}$)$  if $p_{0}=+\infty$ and $p_{1}>1$, then
the embedding of $W$ into $C([0,T_e];B_{1})$ is compact.

\end{Lemma}

In the case that $1\le p_0 <\infty$ and $p_1= 1 $, the lemma is also called the {\bf Generalized Aubin-Lions} which plays a crucial role in the investigation of
the limit as $B\to 0$ in Section~\ref{sec5}.
The proof of  Lemma~\ref{thm2.3}, we refer to, e.g., \cite{8,17,11}. We also will use the lemma on H\"older continuity.
\begin{Lemma} \label{holder-conti}  Let $f (t, x)$ be a function, defined over $Q_{t_e}$, such that\\
 $($i$)$ $f$ is uniformly $($with respect to $x$$)$ H\"older continuous in $t$, with exponent $0 < \alpha\le 1$,
  that is $| f (t, x) - f (s, x)| \le C|t - s|^\alpha$, and\\
 $($ii$)$  $f_x$ is uniformly (with respect to $t$) H\"older continuous in $x$, with exponent $0 < \beta\le 1$,
  that is $| f_x (t, x) - f _x(t, y)| \le C'|y - x|^\beta $.

 Then $f_x$ is uniformly H\"older continuous in $t$ with exponent $\gamma = \alpha \beta/(1 + \beta)$, such that
 $| f_x (t, x) - f _x(s, x)| \le C''|t - s|^\gamma $,  $\forall x \in\bar\Omega,\   0 \le s \le t \le t_e$,
 where $C''$ is a constant which may depend on $C ,\  C'$ and $\alpha, \beta$.

\end{Lemma}

\vskip0.2cm
We now turn back to the proof of Lemma~\ref{lm-est-kappa-2}.

\noindent\textbf{Proof of Lemma~\ref{lm-est-kappa-2}}. Differentiating formally Eq. \eqref{2.1a} (where $\hat\sigma_i^\kappa$
 is replaced by $ \sigma_i^\kappa$) with respect to $t$,  multiplying  by $h_t^{\kappa}$,
using integration by parts, and invoking the
boundary condition \eqref{2.4a}, we obtain
\begin{eqnarray}
 & &\frac{1}{2}\frac{d}{dt}\|h^{\kappa}_t(t)\|^{2}
  + \alpha_{1}\int_{\Omega}\left(  |h_{x}^{\kappa}|_{\kappa}  + B\right)|h _{xt}^{\kappa}|^2    dx \non\\
  &=& - \int_{\Omega} (\alpha_{2}  \sigma_i^\kappa)_t (|  h^\kappa_{x}|_{\kappa}  + B) h^{\kappa}_t
  + (\alpha_{2} \sigma_i^\kappa + \alpha_{3} )(| p|_{\kappa})'|_{p= h^\kappa_{x}}   h^\kappa_{xt} h_t ^{\kappa}dx \non\\
 &=:& I_1+I_2.
 \label{est-kappa-2a}
\end{eqnarray}
Next we are going to estimate $I_1,I_2$.
Noting the periodic boundary conditions, by using the H\"older and Young inequalities
and the Sobolev embedding theorem,   we get
\begin{eqnarray}
 |I_1| &\le& \frac{K\beta}\kappa \left( \int_{\Omega} \left(\int_{\Omega}
   \alpha_{2} |h^\kappa_{xt}(t,x-y)|dy\right)^2\right)^\frac12  \|\,|  h^\kappa_{x}|_{\kappa}  + B\|_{L^\infty(\Omega)} \| h_t ^{\kappa}\|
   \non\\
 &\le & C_\kappa \|h^\kappa_{xt}\| \left(\|  h^\kappa_{xx} + 1\|\,  \| h_t ^{\kappa}\| \right)\non\\
 &\le & \varepsilon \|h^\kappa_{xt}\|^2 + C  (\|  h^\kappa_{xx}\|^2+ 1)  \| h_t ^{\kappa}\|^2.
 \label{est-kappa-2b}
\end{eqnarray}

To evaluate $I_2$, we invoke the Nirenberg inequality in the following form
\begin{eqnarray}
 \| f \|_{L^4} \le C\| f_x \|^\frac14 \| f \|^\frac34  + C'\| f \| ,
 \label{Nirenberg-L4}
\end{eqnarray}
where $f$ will be replaced by $h^\kappa_x$ and $h^\kappa_t$.
It is easy to see that $| (| p|_{\kappa})'|\le 1$, hence applying the Young inequality and recalling
 the estimates \eq{est-kappa-1} we arrive at
\begin{eqnarray}
 |I_2| &\le& C(\|\sigma_i^\kappa\|_{L^4} + 1 ) \| h^\kappa_{xt}\|\,\| h_t ^{\kappa}\|_{L^4}\non\\
 &\le& C(\|h^\kappa_x\|_{L^4} + 1 ) \| h^\kappa_{xt}\|\,\| h_t ^{\kappa}\|_{L^4}\non\\
 &\le& C(\|h^\kappa_{xx}\|^\frac14 \|h^\kappa_x\|^\frac34 + \|h^\kappa_x\| + 1 ) \| h^\kappa_{xt}\|(\|h^\kappa_{tx}\|^\frac14 \|h^\kappa_t\|^\frac34 + \|h^\kappa_t\|)\non\\
 &\le& C(\|h^\kappa_{xx}\|^\frac14   + 1 ) \| h^\kappa_{xt}\|(\|h^\kappa_{tx}\|^\frac14 \|h^\kappa_t\|^\frac34 + \|h^\kappa_t\|)\non\\
 &\le& C\|h^\kappa_{tx}\|^\frac54 ( \|h^\kappa_{xx}\|^\frac14   + 1 ) \|h^\kappa_t\|^\frac34
  + C\|h^\kappa_{tx}\| ( \|h^\kappa_{xx}\|^\frac14   + 1 ) \|h^\kappa_t\|\non\\
 &\le& \varepsilon \|h^\kappa_{tx}\|^2 + C ( \|h^\kappa_{xx}\|^\frac43   + 1 ) \|h^\kappa_t\|^2
  + \varepsilon\|h^\kappa_{tx}\|^2 + C ( \|h^\kappa_{xx}\|^2  + 1 ) \|h^\kappa_t\|^2 .
 \label{est-kappa-2c}
\end{eqnarray}
Combination of \eq{est-kappa-2a}, \eq{est-kappa-2b} and \eq{est-kappa-2c} yields
\begin{eqnarray}
 & &\frac{1}{2}\frac{d}{dt}\|h^{\kappa}_t(t)\|^{2}
  +  \int_{\Omega}\left( \alpha_{1} |h_{x}^{\kappa}|_{\kappa}  + (\alpha_{1} B - 3\varepsilon)\right)|h _{xt}^{\kappa}|^2    dx \non\\
  &\le&   C  \left(\|  h^\kappa_{xx}\|^2 + \|h^\kappa_{xx}\|^\frac43 + 1\right)  \| h_t ^{\kappa}\|^2.
 \label{est-kappa-2d}
\end{eqnarray}

To apply Lemma~\ref{Gronwall}, we define
  $y(t) = \|h^{\kappa}_t(t)\|^{2}$, $A(t)=C  (\|  h^\kappa_{xx}\|^2 + \|h^\kappa_{xx}\|^\frac43 + 1)$ which is
integrable over $[0,t_e]$ by \eq{est-kappa-1}, and $B(t)=0$, then we infer   from \eq{est-kappa-2d} in which we choose
$\varepsilon = \frac{\alpha_{1} B }6 $, that
  \eq{est-kappa-2}, except $\|  h^\kappa_{xx}\|^2\le C$ (which however can be obtained from
  equation \eq{2.1a} with the help of other estimates in \eq{est-kappa-2}), is valid.
  Thus the proof of Lemma~\ref{lm-est-kappa-2} is complete.

\vskip0.2cm
In what follows we will derive more regularities from \eq{est-kappa-1} and \eq{est-kappa-2}.
To this end, we recall  Lemma~\ref{thm2.3},  and let $f=h_x^\kappa$,
$$
 p_{0}=\infty,\ p_{1}=2, \ B_{0}=H^1_{per}(\Omega)\subset\subset B:= C^\alpha_{per}(\bar\Omega),\ B_{1} = L^2(\Omega),
$$
it then follows from  \eq{est-kappa-1} and \eq{est-kappa-2} that
$$
 h_x^\kappa\in L^\infty (0,t_e;B_{0}),\ \frac{\partial h_x^\kappa}{\partial t} =  h^\kappa_{tx} \in L^2(0,t_e;B_{1}),
$$
and   we arrive at
\begin{eqnarray}
 h_x^\kappa\in C([0,t_e];B) = C([0,t_e];C^\alpha_{per}(\bar\Omega)).
 \label{conti-hx}
\end{eqnarray}
Invoking the Sobolev embedding theorem, we also have
\begin{eqnarray}
 |h^\kappa (t,x) - h^\kappa (s,x)|
 &=& \left|\int_s^t h^\kappa_t (\tau,x)d\tau \right|  \le  \left|\int_s^t \|h^\kappa_t(\tau)\| _{L^\infty(\Omega)}d\tau\right| \non\\
 &\le& \left|\int_s^t \|h^\kappa_t(\tau)\| _{H^1_{per}(\Omega)}d\tau\right| \non\\
 &\le& \left(\int_s^t \|h^\kappa_t(\tau)\|^2_{H^1_{per}(\Omega)}d\tau\right)^\frac12
   \left(\int_s^t 1\, d\tau\right)^\frac12 \non\\
 &\le& C | t - s |^\frac12.
 \label{conti-ht}
\end{eqnarray}

\vskip0.2cm
\noindent Completion of the {\bf Proof} of Theorem~\ref{thm2.1}.  Now using \eq{conti-ht}  and \eq{conti-hx} we
may apply  Lemma~\ref{holder-conti}  to conclude that there exists a positive
constant  $0 < \alpha < 1$ such that
$\|h^\kappa_x\|_{ C^{\alpha /2,\alpha} }\le C_\kappa$.
 By the {\it a priori} estimate of the Schauder type for parabolic equations, we thus obtain that
$$
 \|h^\kappa\|_{C^{1+\alpha/2,2+\alpha} ( \bar Q_{ t_e })} \le C_\kappa.
$$
Invoking that $C^{1+\alpha/2,2+\alpha}( \bar Q_{ t_e })\subset\subset C^{\alpha/2,1+\alpha}( \bar Q_{ t_e })$, we
 see the conditions for the Leray-Schauder fixed-point theorem are satisfied.
By definition it is easy to see that $P_0h\equiv 0$.
Thus we are in a position to apply the Leray-Schauder fixed-point theorem (see, e.g., \cite{Gilbarg}),
 and assert that $P_1$ has a fixed point, i.e., $P_1h\equiv h$ and this implies that
 a classical solution exists globally. Hence the proof of Theorem~2.1 is complete.

\section{A prior estimates}\label{subsec2.2}

In this section we are going to derive {\it a prior} estimates for solutions to
the modified problem \eqref{2.1} -- \eqref{2.5}, which are uniform with respect to $\kappa$.
 Since we shall take the limits of approximate solutions as $\kappa\to 0$,
in what follows we may assume that
\begin{eqnarray}
 0<\kappa\le 1,\ \ \gamma H\ {\rm is\  sufficiently\  greater\  than\ }b.
 \label{2.13}
\end{eqnarray}
In this  section, the letter $C$ stands for various positive
constants independent of $\kappa$, but may depend on $B$.

\begin{Lemma} \label{lm2.0} There hold for any $t\in[0,T_{e}]$
\begin{eqnarray}
 \|h^{\kappa}(t)\|^{2} + \int_{0}^{t}\int_{\Omega}\left(\int_{0}^{h_{x}^{\kappa}}|y|_{\kappa} dy + Bh _{x}^{\kappa} \right)
  h_{x}^{\kappa}  dxd\tau
 &\le& C,
 \label{2.16y}\\
  \int_{0}^{t}\|h_{x}^{\kappa} (\tau) \|_{L^3(\Omega)} ^3 d\tau
 &\le& C.
 \label{2.16z}
\end{eqnarray}

\end{Lemma}

\noindent\textbf{Proof.}  Multiplying Eq. \eqref{2.1} by $h^{\kappa}$,
making use of integration by parts, and invoking the
boundary condition \eqref{2.4}, we arrive at
\begin{eqnarray}
 & &\frac{1}{2}\frac{d}{dt}\|h^{\kappa}(t)\|^{2}
  + \alpha_{1}\int_{\Omega}\left(\int_{0}^{h_{x}^{\kappa}}|y|_{\kappa} dy + Bh _{x}^{\kappa} \right)  h_{x}^{\kappa}  dx \non\\
  &=& - \int_{\Omega} (\alpha_{2} \sigma_i^\kappa + \alpha_{3} )(|  h^\kappa_{x}|_{\kappa}  + B) h ^{\kappa}dx
 \non\\
  &=& - \int_{\Omega} \frac{\partial}{\partial x}\int_a^x
  (\alpha_{2} \sigma_i^\kappa + \alpha_{3} )(|  h^\kappa_{x}|_{\kappa}  + B)dy  h ^{\kappa}dx
 \non\\
 &=&  \int_{\Omega}  \int_a^x
  (\alpha_{2} \sigma_i^\kappa + \alpha_{3} )(|  h^\kappa_{x}|_{\kappa}  + B)dy  h_x^{\kappa}dx
 \non\\
 &=:& I.
 \label{2.17z}
\end{eqnarray}
For a term in the left-hand side of \eq{2.17z} we evaluate it as
\begin{eqnarray}
 \int_{\Omega} \int_{0}^{h_{x}^{\kappa}}|y|_{\kappa} dy\, h_{x}^{\kappa} dx
 &\ge& \int_{\Omega} \int_{0}^{h_{x}^{\kappa}}|y| dy\, h_{x}^{\kappa}  dx
 \non\\
 &=& \frac12 \int_{\Omega}  |  h_{x}^{\kappa}| ^{3} dx.
 \label{2.17z1}
\end{eqnarray}
Note that the above inequality is obvious for $h_{x}^{\kappa} \ge 0$, otherwise one may replace $h_{x}^{\kappa}$
by $-|h_{x}^{\kappa}|$ and finds the same inequality.
Applying the Young and H\"older inequalities and Theorem~\ref{thm2.2}, we obtain
\begin{eqnarray}
 |I| & \le &  \int_{\Omega}
  |(\alpha_{2} \sigma_i^\kappa + \alpha_{3} )(|  h^\kappa_{x}|_{\kappa}  + B)|dx \int_{\Omega}  |h_x^{\kappa}|dx
 \non\\
 & \le &   \|\alpha_{2}\sigma_i^\kappa + \alpha_{3}\|_{L^3(\Omega)} \|  \,|h^\kappa_{x}|_{\kappa}  + B\|_{L^3(\Omega)} \|1 \|_{L^3(\Omega)}   \|h_x^{\kappa}\|_{L^3(\Omega)}  \|1 \|_{L^\frac32(\Omega)}
 \non\\
 & \le &  |\Omega| \left(\alpha_{2}\|\sigma_i^\kappa\|_{L^3(\Omega)} + \alpha_{3} |\Omega|\right)
 \left(\|  h^\kappa_{x} \|_{L^3(\Omega)}  + C\right)    \|h_x^{\kappa}\|_{L^3(\Omega)}     .
 \label{2.17y}
\end{eqnarray}
Here $|\Omega|$ denotes the measure of $\Omega$. We now estimate the
term of $\sigma_i^\kappa$. From
estimate \eq{1.6d} it follows that
\begin{eqnarray}
 \|\sigma_{i1}\|_{L^3(\Omega)} \le C \|h_{x} \|_{L^2(\Omega)}
  \le  C \|h_{x} \|_{L^3(\Omega)}.
 \label{2.17y1}
\end{eqnarray}
and from Theorem~\ref{thm2.2} we infer that
\begin{eqnarray}
 \|\sigma_{i2}^\kappa\|_{L^3(\Omega) }\le C_3\|h_x^\kappa\|_{L^3(\Omega)}.
 \label{2.17y2}
\end{eqnarray}
Thus we arrive at
\begin{eqnarray}
 \|\sigma_i^\kappa\|_{L^3(\Omega)}
 & \le &  \|\sigma_{i1}^\kappa\|_{L^3(\Omega)}  + \|\sigma_{i2}^\kappa\|_{L^3(\Omega)}
 \le C \|h_x^\kappa\|_{L^3(\Omega)},
 \label{2.17y3}
\end{eqnarray}
and hence
\begin{eqnarray}
 |I|
 & \le & |\Omega| \left(\alpha_{2}C\|h_x^\kappa\|_{L^3(\Omega)} + \alpha_{3} |\Omega| \right)
 \left(\|  h^\kappa_{x} \|_{L^3(\Omega)}  + C\right)    \|h_x^{\kappa}\|_{L^3(\Omega)}
  \non\\
 & \le & C |\Omega|\alpha_{2}\|h_x^\kappa\|_{L^3(\Omega)} ^{3}
  + \varepsilon  \|h_x^{\kappa}\|^3_{L^3(\Omega)} + C_\varepsilon   .
 \label{2.17y4}
\end{eqnarray}
Combination of inequalities \eq{2.17z} and \eq{2.17y4} yields
\begin{eqnarray}
 & &\frac{1}{2}\frac{d}{dt}\|h^{\kappa}(t)\|^{2}
  + \frac{\alpha_{1}}2\int_{\Omega}{ h_{x}^{\kappa}} \int_{0}^{h_{x}^{\kappa}}|y|_{\kappa} dydx
  + \frac{\alpha_{1}}{ 4}\int_{\Omega} |h_{x}^{\kappa} |^3 dx + B\alpha_{1}\int_{\Omega} |h _{x}^{\kappa} |^2  dx
  \non\\
 & \le &  {\hat\alpha_{2}}\|h_x^\kappa\|_{L^3(\Omega)} ^{3}  + \varepsilon  \|h_x^{\kappa}\|^3_{L^3(\Omega)} + C_\varepsilon ,
 \label{2.17x}
\end{eqnarray}
which implies
\begin{eqnarray}
  \frac{1}{2}\frac{d}{dt}\|h^{\kappa}(t)\|^{2}
  &+& \frac{\alpha_{1}}2\int_{\Omega} { h_{x}^{\kappa}}\int_{0}^{h_{x}^{\kappa}}|y|_{\kappa} dydx
  + \left(\frac{\alpha_{1}}{ 4} -  {\hat \alpha_{2}} - \varepsilon\right) \int_{\Omega} |h_{x}^{\kappa} |^3 dx
  + B\alpha_{1}\int_{\Omega} |h _{x}^{\kappa} |^2  dx \non\\
 & \le &    C_\varepsilon   .
 \label{2.17w}
\end{eqnarray}
Here ${\hat \alpha_{2}} := {\rm meas(\Omega)}C_3\alpha_{2}$.
Therefore choosing that $\alpha_{1}$ is sufficiently greater than $\alpha_{2}$, and $\varepsilon $ is
  suitably small, integrating \eq{2.17w} with respect to $t$, we arrive  at
\begin{eqnarray}
  \|h^{\kappa}(t)\|^{2}
  +   \int_0^t\int_{\Omega} \left(  h_{x}^{\kappa} \int_{0}^{h_{x}^{\kappa}}|y|_{\kappa} dy + |h_{x}^{\kappa} |^3
  +   |h _{x}^{\kappa} |^2 \right) dxd\tau   \le    C  + \|h^{\kappa}_0\|^{2} \le    C .
 \label{2.17v}
\end{eqnarray}
Thus the proof of this lemma is complete.

\begin{Lemma} \label{lm2.1} There holds for any $t\in[0,T_{e}]$
\begin{equation}
 \|h_{x}^{\kappa}(t)\|^{2}+\int_{0}^{t}\int_{\Omega}\left(|h_{x}^{\kappa}|_{\kappa} +B\right) |h_{xx}^{\kappa}|^{2}dxd\tau
   \le C.
 \label{2.16}
\end{equation}

\end{Lemma}

\noindent\textbf{Proof.}  Multiplying Eq. \eqref{2.1} by $-h_{xx}^{\kappa}$,
employing  the technique of  integration by parts with
respect to $x$, and invoking the boundary condition \eqref{2.4},
 we obtain formally for almost all $t$ that
\begin{eqnarray}
 & &\frac{1}{2}\frac{d}{dt}\|h_{x}^{\kappa}(t)\|^{2}
  + \alpha_{1}\int_{\Omega}(|h_{x}^{\kappa}|_{\kappa}
  +B  )|h_{xx}^{\kappa}|^{2}dx
 = \int_{\Omega} (\alpha_{2} \sigma_i^\kappa + \alpha_{3} )(|  h^\kappa_{x}|_{\kappa}  + B) h_{xx}^{\kappa}dx
 \non\\
 &=& \int_{\Omega} \alpha_{2} \sigma_i^\kappa  (|  h^\kappa_{x}|_{\kappa}  + B) h_{xx}^{\kappa} dx
  + \int_{\Omega}\alpha_{3}  (|  h^\kappa_{x}|_{\kappa}  + B) h_{xx}^{\kappa} dx
  \non\\
 &=:& I_1+I_2.
 \label{2.17}
\end{eqnarray}
We may employ the technique of finite difference to justify the formal computation in \eq{2.17}. It is quite
standard so we omit the details.

Now we treat $I_1$ and $I_2$, and first estimate the easier term $I_2$.
Applying the Young inequality with $\varepsilon$, we have
\begin{eqnarray}
 |I_2| &=& \left| \int_{\Omega}\alpha_{3}(|  h^\kappa_{x}|_{\kappa}  + B)^\frac12
 (|  h^\kappa_{x}|_{\kappa}  + B)^\frac12 h_{xx}^{\kappa} dx\right|
  \non\\
 &\le&  C_\varepsilon\int_{\Omega} (|  h^\kappa_{x}|_{\kappa}  + B)dx
  +  \frac\varepsilon2\int_{\Omega} (|  h^\kappa_{x}|_{\kappa}  + B) |   h_{xx}^{\kappa} |^2 dx
  \non\\
 &\le&  C_\varepsilon\int_{\Omega} (|  h^\kappa_{x}| + \kappa + B)dx
  +  \frac\varepsilon2 \int_{\Omega} (|  h^\kappa_{x}|_{\kappa}  + B) |   h_{xx}^{\kappa} |^2 dx
  \non\\
 &\le&  C_\varepsilon\int_{\Omega} (|  h^\kappa_{x}|^2 + C')dx
  +  \frac\varepsilon2 \int_{\Omega} (|  h^\kappa_{x}|_{\kappa}  + B) |   h_{xx}^{\kappa} |^2 dx.
 \label{2.18a}
\end{eqnarray}
Here we used the simple inequality $|p|_{\kappa} \le |p| + \kappa$.  To handle $I_1$, we recall
 Theorem~\ref{thm2.2}, \eq{1.6d}
 then arrive at
\begin{eqnarray}
 |I_1 |
 &=& \alpha_{2} \left|\int_{\Omega} \sigma_i^\kappa  (|  h^\kappa_{x}|_{\kappa}  + B)^\frac12
 \left((|  h^\kappa_{x}|_{\kappa}  + B)^\frac12 h_{xx}^{\kappa} \right)dx\right|
    \non\\
 &\le & \alpha_{2} \left(\int_{\Omega} |\sigma_i^\kappa |^3dx\right)^\frac13
  \left(\int_{\Omega} (|  h^\kappa_{x}|_{\kappa}  + B)^{{\frac12}*6} dx\right)^\frac16
  \left( \int_{\Omega} \left((|  h^\kappa_{x}|_{\kappa}  + B)^\frac12 h_{xx}^{\kappa} \right)^2 dx \right)^\frac12
    \non\\
   &\le & \alpha_{2}\|h^\kappa_{x}\|_{L^3(\Omega)} ^{1+\frac1{ 2}} \left( \int_{\Omega}
   (|  h^\kappa_{x}|_{\kappa}  + B) | h_{xx}^{\kappa} |^2 dx \right)^\frac12
    \non\\
   &\le & C_\varepsilon \|h^\kappa_{x}\|_{L^3(\Omega)} ^{{ 3}} +
    \frac\varepsilon2 \int_{\Omega} (|  h^\kappa_{x}|_{\kappa}  + B) | h_{xx}^{\kappa} |^2 dx .
 \label{2.18b}
\end{eqnarray}
Combining \eq{2.17} with \eq{2.18a} and \eq{2.18b}, integrating it with respect to $t$, and
making use of  the Young inequality, we then arrive at
\begin{eqnarray}
 &&\frac{1}{2} \|h_{x}^{\kappa}(t)\|_{2}^{2}
  + \alpha_{1}\int_0^t\int_{\Omega}(|h_{x}^{\kappa}|_{\kappa} + B) |h_{xx}^{\kappa}|^{2}dxd\tau
   \nonumber\\
 &\le & C_\varepsilon \int_0^t\|h^\kappa_{x}\|_{L^3(\Omega)} ^{3} d\tau + C +
     \varepsilon \int_0^t \int_{\Omega} (|  h^\kappa_{x}|_{\kappa}  + B) | h_{xx}^{\kappa} |^2 dx d\tau.
  \label{2.18}
\end{eqnarray}
Now we choose $\varepsilon$ small enough so  that $\alpha_{1}- \varepsilon>0$, recall the estimates in Lemma~\ref{lm2.0}, then
\begin{eqnarray}
 &&\frac{1}{2} \|h_{x}^{\kappa}(t)\|_{2}^{2}
  + (\alpha_{1} - \varepsilon)\int_0^t\int_{\Omega}(|h_{x}^{\kappa}|_{\kappa} + B) |h_{xx}^{\kappa}|^{2}dxd\tau
   \nonumber\\
 &\le & C_\varepsilon \int_0^t\|h^\kappa_{x}\|_{L^3(\Omega)} ^{3} d\tau + C \nonumber\\
 &\le & C .
 \label{2.19}
\end{eqnarray}
Therefore, the proof of Lemma~\ref{lm2.1} is complete.

\begin{Corollary}\label{cor2.1} There hold for any $t\in[0,T_e]$
\begin{eqnarray}
 \int_{0}^{t}\int_{\Omega}(|h_{x}^{\kappa}|_{\kappa}|h_{xx}^{\kappa}|)^{\frac{4}{3}}dxd\tau
 &\le& C,
 \label{2.22a}\\
 \int_{0}^{t}\int_{\Omega}(|h_{x}^{\kappa} h_{xx}^{\kappa}|)^{\frac{4}{3}}dxd\tau
 &\le& C,
 \label{2.22b}\\
 \int_0^t\| ( h_{x}^{\kappa} )^2  \|^\frac43 _{W^{1,\frac{4}{3} } (\Omega) }d\tau
 &\le& C,
 \label{2.22c}\\
 \int_0^t\|  h_{x}^{\kappa}  \|^\frac83 _{L^{\infty} (\Omega) }d\tau
 &\le& C.
 \label{2.22}
\end{eqnarray}

\end{Corollary}

\noindent\textbf{Proof.}  By H\"{o}lder's inequality, we have for some $1\le p<2, q=\frac{2}{p}$ and $\frac{1}{q}+\frac{1}{q'}=1$ that
\begin{eqnarray}
  & & \int_{0}^{t}\int_{\Omega}(|h_{x}^{\kappa}|_{\kappa}|h_{xx}^{\kappa}|)^{p}dxd\tau
 = \int_{0}^{t}\int_{\Omega}|h_{x}^{\kappa}|_{\kappa}^{\frac{p}{2}}
 \left(|h_{x}^{\kappa}|_{\kappa}^{\frac{p}{2}}|h_{xx}^{\kappa}|^{p}\right)dxd\tau
 \nonumber\\
 &\le&\Big(\int_{0}^{t}\int_{\Omega}|h_{x}^{\kappa}|_{\kappa}^{\frac{pq'}{2}}dxd\tau\Big)^{\frac{1}{q'}}
 \Big(\int_{0}^{t}\int_{\Omega}|h_{x}^{\kappa}|_{\kappa}^{\frac{pq}{2}}|h_{xx}^{\kappa}|^{pq}dxd\tau\Big)^{\frac{1}{q}}
 \nonumber\\
 &\le&\Big(\int_{0}^{t}\int_{\Omega}|h_{x}^{\kappa}|_{\kappa}^{\frac{p}{2-p}}dxd\tau\Big)^{\frac{2-p}{2}}
 \Big(\int_{0}^{t}\int_{\Omega}|h_{x}^{\kappa}|_{\kappa}|h_{xx}^{\kappa}|^{2}dxd\tau\Big)^{\frac{p}{2}}.
 \label{2.23}
\end{eqnarray}
Inequality \eqref{2.16} implies for $\frac{p}{2-p}\leq2$, i.e., $p\le \frac{4}{3}$,
that the right-hand side of \eqref{2.23} is bounded, hence \eq{2.22a} is true.

Invoking the basic fact that $|p|_\kappa\ge |p|$, from \eq{2.22a} it follows that \eq{2.22b} holds.
To prove \eq{2.22}, we apply the Poincar\'e inequality in the following form
$$
 \|f-\bar f\|_{L^p(\Omega)} \le C\|f_x\|_{L^p(\Omega)},
$$
where $\bar f = \int_\Omega f(x)dx/|\Omega|$, choose $p=\frac43$,
then   from \eq{2.22b} we deduce that
\begin{eqnarray}
 \|(h_{x}^{\kappa})^2 - \overline{ (h_{x}^{\kappa})^2 }\|_{W^{1,\frac43}(\Omega)}
 &\le& \|\left((h_{x}^{\kappa})^2 - \overline{ (h_{x}^{\kappa})^2 }\right)_x\|_{W^{1,\frac43}(\Omega)}
 = \|\left((h_{x}^{\kappa})^2 \right)_x\|_{W^{1,\frac43}(\Omega)}   \non\\
 &=& 2\|  h_{x}^{\kappa} h_{xx}^{\kappa}\|_{W^{1,\frac43}(\Omega)} ,
 \label{2.23a}
\end{eqnarray}
hence
\begin{eqnarray}
 \int_0^t \|(h_{x}^{\kappa})^2 - \overline{ (h_{x}^{\kappa})^2 }\|^\frac43_{W^{1,\frac43}(\Omega)} d\tau
 \le   2\int_0^t\|  h_{x}^{\kappa} h_{xx}^{\kappa}\|^\frac43_{W^{1,\frac43}(\Omega)} d\tau\le C,
 \label{2.23b}
\end{eqnarray}
which implies
\begin{eqnarray}
 \int_0^t \|(h_{x}^{\kappa})^2  \|^\frac43_{W^{1,\frac43}(\Omega)} d\tau
 &\le& \int_0^t \|(h_{x}^{\kappa})^2 - \overline{ (h_{x}^{\kappa})^2 }\|^\frac43_{W^{1,\frac43}(\Omega)} d\tau
  + \int_0^t \|  \overline{ (h_{x}^{\kappa})^2 }\|^\frac43_{W^{1,\frac43}(\Omega)} d\tau
  \non\\
 &\le& C + \int_0^t\overline{ (h_{x}^{\kappa})^2 }\| 1 \|^\frac43_{L^{\frac43}(\Omega)} d\tau
 \non\\
 &\le& C + \sup_{0\le \tau\le t} \overline{ (h_{x}^{\kappa})^2 }(\tau) \int_0^t\| 1 \|^\frac43_{L^{\frac43}(\Omega)} d\tau
 \non\\
 &\le& C.
 \label{2.23c}
\end{eqnarray}
By the Sobolev embedding theorem we have $W^{1,\frac43}(\Omega)\subset L^{\infty}(\Omega)$,
and \eq{2.22c} follows.
Hence the proof of the corollary is complete.

\begin{Lemma} \label{lm2.2} There holds for any $t\in[0,T_{e}]$
\begin{equation}
 \|h_{t}^{\kappa}\|_{L^{\frac{4}{3}}(Q_{T_{e}})}\le C.
 \label{2.24}
\end{equation}

\end{Lemma}

\noindent\textbf{Proof.} Recalling the regularity about $h_{t}^\kappa$  we use the integration by
 parts and obtain
\begin{eqnarray}
   (h^\kappa,\varphi_{t})_{Q_{T_{e}}}
   &=& \int_0^{T_e}\frac{d}{dt}(h^\kappa,\varphi )_{\Omega}dt - (h_{t}^\kappa,\varphi)_{Q_{T_{e}}}
   \non\\
   &=& \left.(h^\kappa,\varphi )_{\Omega}\right|_0^{T_e} - (h_{t}^\kappa,\varphi)_{Q_{T_{e}}}
  \non\\
   &=& - (h^\kappa(0),\varphi (0))_{\Omega}  - (h_{t}^\kappa,\varphi)_{Q_{T_{e}}}
 \label{1.11c}
\end{eqnarray}
thus one has
\begin{eqnarray}
   (h_{t}^\kappa,\varphi)_{Q_{T_{e}}}
   = \Big(\big({\alpha_{1}} h^\kappa _{xx} - \left(\alpha_{2} \sigma_i^\kappa + \alpha_{3}
    \right)\big) (|h^\kappa _{x}|_\kappa + B), \varphi\Big) _{Q_{T_{e}}}  .
 \label{1.11d}
\end{eqnarray}
Making use of Theorem~\ref{thm2.2} with $p=\frac43$, and estimates
\eq{2.22a}, \eq{2.22},  and the H\"older inequality, we  have
\begin{eqnarray}
  && | (h_{t}^\kappa,\varphi)_{Q_{T_{e}}}|\non\\
  &\le& C\Big(\| \, |h^\kappa _{x}|_\kappa h^\kappa _{xx}\|_{L^\frac43(Q_{T_e})}
  + \| h^\kappa _{xx}\|_{L^\frac43(Q_{T_e})}  \Big)\| \varphi \|_{L^4(Q_{T_e})}
     \non\\
  & & + C\Big(  \|  \sigma_i^\kappa\|_{L^\frac83(0,{T_e},L^\frac43(\Omega))}
     \|  h^\kappa _{x}\|_{L^\frac83(0,{T_e};L^\infty(\Omega))}
      +1 \Big)\| \varphi \|_{L^4(Q_{T_e})} \non\\
  &\le&   C\Big(  1  + \|  h^\kappa _{x} \|_{L^\frac83(0,{T_e},L^{ \infty}(\Omega))}  \Big)\| \varphi \|_{L^4(Q_{T_e})}\non\\
  &\le&   C \| \varphi \|_{L^4(Q_{T_e})}
 \label{1.11e}
\end{eqnarray}
for all $\varphi \in {L^4(Q_{T_e})} $. Thus \eq{1.11e} implies
   that  $h_{t}^{\kappa}\in L^{\frac{4}{3}}(Q_{T_{e}})$
and   \eqref{2.24} holds. The proof of the lemma is complete.

\section{Existence of solutions to the original IBVP}

  In this  section we shall use the {\it a prior}
  estimates established in Section~2.2 to
  investigate the convergence of $h^{\kappa}$ as $\kappa\to 0$, and
   show that
  there exists a subsequence, which converges to a weak solution to  the
   initial-boundary value problem \eqref{1.1} -- \eqref{1.3}, thereby
   prove Theorem~\ref{thm1.1}.

\begin{Lemma}\label{lm2.3} There exists a subsequence of $h_{x}^{\kappa}$
(we still denote it by $h_{x}^{\kappa}$) such that
\begin{eqnarray}
 h_{x}^{\kappa}
 &\to& h_{x} {\rm \ strongly\ in\ } L^{2}(Q_{T_{e}}),
 \label{2.25}\\
 |h_{x}^{\kappa}|_{\kappa}
 &\to& |h_{x}| {\rm \ strongly\ in\ }  L^{2}(Q_{T_{e}}),
 \label{2.26}\\
 ~~|h_{x}^{\kappa}|_{\kappa}h_{x}^{\kappa}
 &\to& |h_{x}|h_{x} {\rm \ strongly\ in\ }  L^{1}(Q_{T_{e}})
 \label{2.27}
\end{eqnarray}
as $\kappa\to 0$.

\end{Lemma}

\noindent\textbf{Proof.}  Let $p_{0}=2, p_{1}=\frac{4}{3}$ and
$$
 B_{0}=H^{1}_{{\rm per}}(\Omega), \ B_{1}=L^{2}(\Omega), \ B_{2}=W^{-1,\frac{4}{3}}_{{\rm per}}(\Omega).
$$
These spaces satisfy the assumptions of Lemma~\ref{thm2.3}. Since
estimate  \eqref{2.16} implies that
$ h_{xx}^{\kappa}\in L^{2}(0,T_{e};L^{2}(\Omega))$, then
\begin{equation}
 h_{x}^{\kappa}\in L^{2}(0,T_{e};H^{1}_{{\rm per}}(\Omega)).
 \label{2.28}
\end{equation}
From estimate \eqref{2.24}, we have
\begin{equation}
 h_{xt}^{\kappa}\in L^{\frac{4}{3}}(0,T_{e};W_{{\rm per}}^{-1,\frac{4}{3}}).
 \label{2.29}
\end{equation}
It follows from Theorem~\ref{thm2.3} that
$$
 h_{x}^{\kappa}\to h_{x}\ {\rm  strongly\ in}\ L^{2}(Q_{T_{e}}),
$$
 as $\kappa\to 0$.
  This proves \eqref{2.25}.

It is easy to see that
$$
 |\sqrt{x} - \sqrt{y}| \le \sqrt{|x-y|}
$$
for all $x,y\in\cR^+$. From this we deduce that
as $ \kappa\to 0 $
\begin{eqnarray}
 \big\||h_{x}^{\kappa}|_{\kappa}-|h_{x}|\big\|^2_{L^{2}(Q_{T_{e}})}
 &\le& \big\| \sqrt{\left|(h_{x}^{\kappa})^2 +\kappa^2 - (h_{x})^2\right|}\big\|^2_{L^{2}(Q_{T_{e}})}
 \non\\
 &\le&  \int_{Q_{T_{e}} }\left(\left|(h_{x}^{\kappa})^2 -( h_{x})^2\right| + \kappa^{2} \right)dxd\tau
 \non\\
 &\le&  \|h_{x}^{\kappa} + h_{x} \|_{L^{2}(Q_{T_{e}})}  \|h_{x}^{\kappa} - h_{x} \|_{L^{2}(Q_{T_{e}})}
 + \big\|\kappa^2\big\|_{L^{2}(Q_{T_{e}})}
 \non\\
 &\le&  C  \|h_{x}^{\kappa} - h_{x} \|_{L^{2}(Q_{T_{e}})}
 + \big\|\kappa^2\big\|_{L^{2}(Q_{T_{e}})}
 \non\\
 &\to& 0.
\end{eqnarray}
 From this
we infer that $|h_{x}^{\kappa}|_{\kappa}$ converges to
$|h_{x}|$ strongly  in $L^{2}(Q_{T_{e}})$ as $\kappa\to 0$.
This proves \eqref{2.26}. Combining \eqref{2.25} and \eqref{2.26},
we get \eqref{2.27} immediately.

\vskip0.32cm
\noindent{\textbf{Proof of Theorem~\ref{thm1.1}.}} We have
$\|h^{ \kappa}\|_{L^{\infty} (0,T_{e};H_{{\rm per}}^{1}(\Omega))}\le  C$,
and $\|h^{ \kappa}\|_{L^{2}(0,T_{e};H^{2}_{{\rm per}}(\Omega))}  \le  C$
by \eqref{2.16}. This implies
$h\in L^{\infty}(0,T_{e};H_{{\rm per}}^{1}(\Omega))\cap L^{2}(0,T_{e};H^{2}_{{\rm per}}(\Omega))$,
since we can select a subsequence of $h^{ \kappa}$ which converges
weakly to $h$ in this space. Thus, $h$ satisfies \eqref{1.10}.

 It therefore suffices to show that problem \eqref{1.1} -- \eqref{1.3}
is fulfilled in the weak sense which means we need to prove
 the relation \eqref{1.11} holds. To
 this end, we employ the following equality
\begin{eqnarray}
  (h^\kappa,\varphi_{t})_{Q_{T_{e}}} -  {\alpha_{1}} \left( \int_{0}^{h^\kappa _{x}}|p |_\kappa dp
   + B h^\kappa_{x},\varphi_{x} \right)_{Q_{T_{e}}}
   &-&  \left((\alpha_{2} \sigma_i^\kappa + \alpha_{3} )(|h^\kappa _{x}|_\kappa + B), \varphi\right) _{Q_{T_{e}}}
 \nonumber\\
 &+& (h_{0},\varphi(0))_{\Omega}=0.
\end{eqnarray}
From which we see that equation \eqref{1.11} follows   if we show that
\begin{eqnarray}
 (h^{\kappa},\varphi_{t})_{Q_{T_{e}}} &\to& (h,\varphi_{t})_{Q_{T_{e}}},
 \label{2.32}\\
 \Big(\int_{0}^{h_{x}^{\kappa}}|y|_{\kappa},\varphi_{x}\Big)_{Q_{T_{e}}} &\to&
 \Big(\frac{1}{2}|h_{x}|h_{x},\varphi_{x}\Big)_{Q_{T_{e}}},
 \label{2.33}\\
 (|h_{x}^{\kappa}|_{\kappa},\varphi)_{Q_{T_{e}}} &\to& (|h_{x}|,\varphi)_{Q_{T_{e}}},
 \label{2.34}\\
 (h_{x}^{\kappa},\varphi_{x})_{Q_{T_{e}}} &\to& (h_{x},\varphi_{x})_{Q_{T_{e}}},
 \label{2.35}\\
 (\sigma_{i}^{\kappa}|h_{x}^{\kappa}|_{\kappa},\varphi)_{Q_{T_{e}}}
  &\to& (\sigma_{i}|h_{x}|,\varphi)_{Q_{T_{e}}}
 \label{2.36}
\end{eqnarray}
as $\kappa\to0$.
Now, the conclusions \eqref{2.32} and \eqref{2.35}
follow easily from \eqref{2.16}, and
the relation \eqref{2.34} follows
from \eqref{2.26}.
It remains to prove
\eqref{2.33} and \eqref{2.36}. To prove \eqref{2.33} we write
\begin{eqnarray}
 \int_{0}^{h_{x}^{\kappa}}|y|_{\kappa}dy-\frac{1}{2}|h_{x}|h_{x}
 &=&\Big(\int_{0}^{h_{x}^{\kappa}}|y|_{\kappa}dy-\frac{1}{2}|h_{x}^{\kappa}|_{\kappa}h_{x}^{\kappa}\Big)
 +\frac{1}{2}(|h_{x}^{\kappa}|_{\kappa}h_{x}^{\kappa}-|h_{x}|h_{x})
 \nonumber\\
 &:=&I_{1}+I_{2}.
 \label{2.37}
\end{eqnarray}
The conclusion \eqref{2.27} implies
\begin{eqnarray}
 \|I_{2}\|_{L^{1}(Q_{T_{e}})}\to 0,
 \label{2.38}
\end{eqnarray}
as $\kappa\to0$. Next we handle $I_{1}$ as follows.
\begin{eqnarray}
 |I_{1}| &=& \Big|\int_{0}^{h_{x}^{\kappa}}|y|_{\kappa}dy-\frac{1}{2}|h_{x}^{\kappa}|_{\kappa}h_{x}^{\kappa}\Big|
  \le  \Big|\int_{0}^{h_{x}^{\kappa}}|y|_{\kappa}dy-\int_{0}^{h_{x}^{\kappa}}|y|dy\Big|
 \nonumber\\
 &\le&  \int_{0}^{|h_{x}^{\kappa}|}\big||y|_{\kappa}-|y|\big|dy
 \nonumber\\
 &\le& \int_{0}^{|h_{x}^{\kappa}|}\kappa dy=\kappa|h_{x}^{\kappa}|,
\end{eqnarray}
whence \eqref{2.16} implies
\begin{eqnarray}
 \|I_{1}\|_{L^{1}(Q_{T_{e}})} &\le& C\|I_{1}\|_{L^{\infty}(0,T_{e};L^{2}(\Omega))}
 \nonumber\\
 &\le& C\|I_{1}\|_{L^{\infty}(0,T_{e};H^{1}(\Omega))} \le  C\kappa
  \to 0,
\end{eqnarray}
as $\kappa\to0$. From this relation
and \eqref{2.37}, \eqref{2.38} we obtain
\begin{eqnarray}
 \Big\|\int_{0}^{h_{x}^{\kappa}}|y|_{\kappa}dy-\frac{1}{2}|h_{x}|h_{x}\Big
 \|_{L^{1}(Q_{T_{e}})}\to 0,
\end{eqnarray}
which implies \eqref{2.33}.
Finally we prove \eqref{2.36}. Applying the compactness lemma and Theorem~\ref{thm2.2} with $p=2$
we get that
\begin{eqnarray}
 \sigma_{i}^{\kappa}\to\sigma_{i}\ {\rm strongly\ in}\ L^{2}(Q_{T_{e}})
\end{eqnarray}
where
$$
 \sigma_{i}(t,x) = P.V.\int_{-\infty}^\infty\frac{K\beta h_x(t,x_1)}{ x - x_1  }dx_1 .
$$
Then recalling \eq{2.26} one concludes that
\begin{eqnarray}
 \sigma_{i}^{\kappa}|h_{x}^{\kappa}|_{\kappa}
 \to\sigma_{i}|h_{x}|\ {\rm strongly\ in}\ L^{1}(Q_{T_{e}}),
\end{eqnarray}
which implies \eqref{2.36}. Thus \eqref{1.11} holds.

We now investigate the regularity
properties of the solution  stated in \eqref{1.12a} and \eqref{1.12}. For  \eqref{1.12a},
we apply the estimates in Lemmas 3.1,  3.2 and the relation \eq{2.28}.
 The assertion  $h_{t}\in L^{\frac{4}{3}}(Q_{T_{e}})$
is implied by \eqref{2.24}. To verify the second assertion in \eqref{1.12},
we use estimates \eq{2.22b} and \eq{2.22} in Corollary 3.1, and also the strong convergence
 result in Lemma 4.1.
Consequently \eqref{1.12} holds.

It remains to prove the uniqueness. To this end, we recall the
regularity of $h^\kappa _t$, and definition  \eq{1.11}, using
integration by parts, to get
\begin{eqnarray}
  &&-(h_{t},\varphi)_{Q_{T_{e}}} -
   {\alpha_{1}} \left( \frac12h _{x} | h _{x} |
   + B h_{x},\varphi_{x} \right)_{Q_{T_{e}}} - \left((\alpha_{2}
   \sigma_i + \alpha_{3} )(|h_{x}| + B), \varphi\right) _{Q_{T_{e}}}
 \nonumber\\
 & = & -(h(0),\varphi{ (0)} )_{ {\Omega}} + (h_{0},\varphi(0))_{\Omega} =0.
 \label{1.11z}
\end{eqnarray}
Suppose that there exist two solutions $h_1,h_2$, let $h = h_1 - h_2$,
then from \eq{1.11z} we infer that
\begin{eqnarray}
  && (h_{t},\varphi)_{Q_{T_{e}}} +
   \frac{\alpha_{1}}2 \left(h _{1 x} | h _{1x} | - h _{2x} | h _{2x} |,\varphi_{x} \right)_{Q_{T_{e}}}
   + {\alpha_{1}} \left(B h_{x},\varphi_{x} \right)_{Q_{T_{e}}}
  \nonumber\\
  & &  + \left((\alpha_{2}
   \sigma_i^1 + \alpha_{3} )(|h_{1x}| + B) - (\alpha_{2}
   \sigma_i^2 + \alpha_{3} )(|h_{2x}| + B), \varphi\right) _{Q_{T_{e}}}
 \nonumber\\
 & = &  0.
 \label{1.11y}
\end{eqnarray}
Here $\sigma_i^j$ with $j=1,2$ stand for  the
formulas of $\sigma_i$ in which $h$   is replaced by $h^j$, respectively.
Since $C_0^\infty(Q_t$) is dense in $L^2(Q_t)$, we can choose $\varphi = h$,
using the monotonicity property
$$
 (x|x| - y|y|)(x-y) \ge 0
$$
to infer from  \eq{1.11y} that
\begin{eqnarray}
  &&  \frac12 \|h(t)\|^2
   + {\alpha_{1}} B\| h_{x} \|^2_{Q_{T_{e}}}
   + \left((\alpha_{2}
   \sigma_i^1 + \alpha_{3} )(|h_{1x}| + B) - (\alpha_{2}
   \sigma_i^2 + \alpha_{3} )(|h_{2x}| + B), h\right) _{Q_{T_{e}}}
 \nonumber\\
 & \le &  \frac12 \|h(0)\|^2 =0 .
 \label{1.11x}
\end{eqnarray}
We write
\begin{eqnarray}
  &I:= & \left((\alpha_{2}\sigma_i^1 + \alpha_{3} )(|h_{1x}| + B) - (\alpha_{2}
   \sigma_i^2 + \alpha_{3} )(|h_{2x}| + B), h\right) _{Q_{T_{e}}}
 \nonumber\\
 &=& \left(\int_a^x(\alpha_{2}\sigma_i^1 + \alpha_{3} )(|h_{1x}| - |h_{2x}|)
  + \alpha_{2}(\sigma_i^1 - \sigma_i^2  )(|h_{2x}| + B)dy, h_x \right) _{Q_{T_{e}}} ,
 \end{eqnarray}
 whence applying again Theorem~\ref{thm2.2} and the H\"older inequality,   we obtain
 \begin{eqnarray}
 |I|& \le & C\int_0^t  \int_a^d\left(( |\sigma_i^1| + 1 ) |h_{x}|
  + |\sigma_i^1 - \sigma_i^2 | (|h_{2x}| + B)\right)dy\,  \|h_x \|  d\tau
 \nonumber\\
 & \le & C\int_0^t \left( (\|\sigma_i^1\| + 1 ) \|h_{x}\|
  + \|\sigma_i^1 - \sigma_i^2 \| (\|h_{2x}\| + B) \right)  \|h_x \|  d\tau
 \nonumber\\
 & \le & C\int_0^t \left(  (\|h_{1x}\| + \|h_{2x}\| + 1) \right)  \|h_x \|^2  d\tau
 \nonumber\\
 & \le & C'\int_0^t   \|h_x \|^2  d\tau.
 \label{1.11w}
\end{eqnarray}

Now  choosing $\alpha_{1}$ sufficiently large, we infer from \eq{1.11w} and
 \eq{1.11x} that
\begin{eqnarray}
   \frac12 \|h(t)\|^2
   + ({\alpha_{1}} B -C')\| h_{x} \|^2_{Q_{T_{e}}} \le 0 ,
 \label{1.11u}
\end{eqnarray}
hence  $\|h(t)\|=0$ which implies that $h=0$  for almost
all $(t,x)\in Q_{T_{e}}$, the uniqueness follows, and
  thus the proof of Theorem~\ref{thm1.1} is complete.

\section{The limit of $h_B$ as $B$ vanishes}
\label{sec5}

This section  is devoted to the investigation the limit of $h_B$ as $B\to 0$,
and to the proof of Theorem~\ref{thm1.1B}.  We will denote the solution $h$ to problem \eq{1.1} -- \eq{1.3} by $h_B$.
Thus we need {\it a priori} estimates which are independent of $B$ and $B$ may be taken to meet
$$
 0<B\le 1.
$$
Those estimates in Lemmas~\ref{lm2.0} and \ref{lm2.1},  and
 Corollary~\ref{cor2.1} are of this type. In this section a universal constant $C$
  is independent of $B$.

To prove Theorem~\ref{thm1.1B}, we shall obtain more estimates
as follows.
\begin{Lemma} \label{lm-est-B} There hold for any $t\in [0,T_e]$ and
for any $\phi\in L^\infty(0,T_e;H^2_{\rm per}(\Omega))$ that
\begin{eqnarray}
 |\left((|h^\kappa_{x}|h^\kappa_{x} )_t,\phi\right)|
 & \le& C\|\phi\|_{L^\infty(0,T_e;H^2_{\rm per}(\Omega))},
  \label{est-B1}\\
 \|(|h^\kappa_{x}|h^\kappa_{x} )_t\|_{L^1(0,t;H^{-2}_{\rm per}(\Omega))}
 & \le& C.
  \label{est-B}
\end{eqnarray}

\end{Lemma}

\vskip0.32cm
\noindent{\textbf{Proof.}}  For a rigorous procedure, we derive   estimate \eq{est-B1} from Eq. \eq{2.1a}, where $\hat\sigma^\kappa_i$  is replaced by
$\sigma^\kappa_i$), thus we see
the solution $h$ depends on both $\kappa$ and $B$. We
shall write $h=h^\kappa_B(t,x)$. However
as in Section~4 one can pass $h^\kappa_B$  to its limit as $\kappa\to 0$, and
get solutions $h_B$, and hence we get estimates for $h_B$.
 To this end, we take an arbitrary $\phi\in L^\infty(0,t;H^2_{\rm per}(\Omega))$,
 multiply $h^\kappa_t$ by $(|h^\kappa_x|_\delta\phi)_x$, and
 integrate the production with respect to $x,t$, and arrive at
\begin{eqnarray}
 I &:=& 2\int_0^t\int_\Omega h^\kappa_t(|h^\kappa_x|_\delta\phi)_x  dxd\tau
   = 2\int_0^t\int_\Omega h^\kappa_t\left(\frac{h^\kappa_x h^\kappa_{xx}}{|h^\kappa_x|_\delta} \phi +  |h^\kappa_x|_\delta \phi_x\right) dxd\tau \non\\
   &:=& I_1+I_2.
  \label{est-B2}
\end{eqnarray}
Here $|h^\kappa_x|_\delta = \sqrt{|h^\kappa_x|^2+\delta^2}$ with a small positive parameter $\delta\le 1$. For the sake of notations' simplicity, we still denote
$h=h^\kappa_B$ by $h=h^\kappa$.

We now treat $I_1$ and $I_2$. First for $I_1$ we invoke Eq. \eq{2.1a} (where $\hat\sigma^\kappa_i$  is replaced by
$\sigma^\kappa_i$),
to get
\begin{eqnarray}
 |I_1|
 &\le&  2|\int_0^t\int_\Omega h^\kappa_t \frac{h^\kappa_x h^\kappa_{xx}}{|h^\kappa_x|_\delta} \phi dx |
 \le   C\int_0^t\int_\Omega |h^\kappa_t   h^\kappa_{xx}  \phi |dxd\tau  \non\\
 &\le&   C\int_0^t\int_\Omega\left(( |h^\kappa_x|_\kappa + B)   |h^\kappa_{xx}|^2|  \phi | + (|\sigma^\kappa_i| + 1)( |h^\kappa_x|_\kappa + B)   |h^\kappa_{xx}   \phi |\right)dxd\tau  \non\\
 &:=& I_{11} + I_{12} .
  \label{est-B2}
\end{eqnarray}
By using the estimates in Lemmas~\ref{lm2.0} and \ref{lm2.1},  and
 Corollary~\ref{cor2.1}, one gets easily
\begin{eqnarray}
 |I_{11}|
 &\le&  C\|  \phi \|_{L^\infty(Q_t)}\int_0^t\int_\Omega ( |h^\kappa_x|_\kappa + B)   |h^\kappa_{xx}|^2 dxd\tau  \non\\
 &\le& C \|  \phi \|_{L^\infty(Q_t)}
 \le  C \|  \phi \|_{L^\infty(0,t;H^2_{\rm per}(\Omega))},
  \label{est-B2a}
\end{eqnarray}
and
\begin{eqnarray}
 |I_{12}|
 &\le& C\|  \phi \|_{L^\infty(Q_t)}\int_0^t \left((\|\sigma^\kappa_i\|_{L^3(\Omega)} + 1)\|\, |h^\kappa_x|^\frac12_\kappa + B\|_{L^6(\Omega)} \|  (|h^\kappa_x|^\frac12_\kappa + B)|h^\kappa_{xx}\|   \right) d\tau  \non\\
 &\le& C \|  \phi \|_{L^\infty(0,t;H^2_{\rm per}(\Omega))}
  \int_0^t \left(\left( \|  h^\kappa_x\|^\frac32_{L^3(\Omega)} + 1 \right) \|  (|h^\kappa_x|^\frac12_\kappa + B)|h^\kappa_{xx}\|   \right) d\tau  \non\\
 &\le& C \|  \phi \|_{L^\infty(0,t;H^2_{\rm per}(\Omega))} \int_0^t \left( \left( \|  h^\kappa_x\|^\frac32_{L^3(\Omega)} + 1 \right)^2 +  \|  (|h^\kappa_x|^\frac12_\kappa + B)|h^\kappa_{xx}\|^2  \right)d\tau\non\\
 &\le& C \|  \phi \|_{L^\infty(0,t;H^2_{\rm per}(\Omega))}.
  \label{est-B2b}
\end{eqnarray}

Next $I_2$ is evaluated as follows.
\begin{eqnarray}
 |I_{2} |
 &\le&  C\int_0^t\int_\Omega  |h^\kappa_{t}|  |h^\kappa_x|_\delta    | \phi_{x} | dxd\tau  \non\\
 &\le&  C\int_0^t\int_\Omega  \left(( |h^\kappa_x|_\kappa + B)   |h^\kappa_{xx}| + (|\sigma^\kappa_i| + 1)( |h^\kappa_x|_\kappa + B)    \right)  |h^\kappa_x|_\delta    | \phi_{x} | dxd\tau\non\\
 &\le&  C\| \phi_{x} \|_{L^\infty(Q_t)} \int_0^t\int_\Omega  \Big(( |h^\kappa_x|_\kappa + B)   |h^\kappa_{xx}| + (|\sigma^\kappa_i| + 1)( |h^\kappa_x|_\kappa + B)    \Big)  (|h^\kappa_x|+1)    dxd\tau \non\\
 &:=& I_{21} + I_{22}.
  \label{est-B3}
\end{eqnarray}
Using the estimates in Lemmas~\ref{lm2.0} and \ref{lm2.1}, we obtain
\begin{eqnarray}
 |I_{21} |
 &\le&  C\| \phi_{x} \|_{L^\infty(Q_t)} \int_0^t\int_\Omega  ( |h^\kappa_x|_\kappa + B)^\frac12  |h^\kappa_{xx}|    (|h^\kappa_x|+1)^\frac32    dxd\tau \non\\
 &\le&  C\| \phi_{x} \|_{L^\infty(Q_t)} \int_0^t\int_\Omega \left(  ( |h^\kappa_x|_\kappa + B)   |h^\kappa_{xx}|^2 +     (|h^\kappa_x|+1)^3  \right)  dxd\tau \non\\
 &\le&  C\| \phi\|_{L^\infty(0,t;H^2_{\rm per}(\Omega))},
  \label{est-B3a}
\end{eqnarray}
and  from the estimates in
 Corollary~\ref{cor2.1} it follows that
 \begin{eqnarray}
 |I_{22} |
 &\le&  C\| \phi_{x} \|_{L^\infty(Q_t)} \int_0^t\int_\Omega    (|\sigma^\kappa_i| + 1)  (|h^\kappa_x|^2 + 1)    dxd\tau \non\\
 &\le&  C\| \phi\|_{L^\infty(0,t;H^2_{\rm per}(\Omega))} \int_0^t\int_\Omega    (|\sigma^\kappa_i| + 1)^3 +  (|h^\kappa_x|^2 + 1)^\frac32     dxd\tau \non\\
 &\le&  C\| \phi\|_{L^\infty(0,t;H^2_{\rm per}(\Omega))}.
  \label{est-B3b}
\end{eqnarray}

Therefore combining \eq{est-B2} -- \eq{est-B3b} together we arrive at
\begin{eqnarray}
 |I | \le  C\| \phi\|_{L^\infty(0,t;H^2_{\rm per}(\Omega))}.
  \label{est-B4}
\end{eqnarray}

We now rewrite $I$ as
\begin{eqnarray}
 I &:=& 2\int_0^t\int_\Omega h^\kappa_t(|h^\kappa_x|_\delta\phi)_x dxd\tau
   = - 2\int_0^t\int_\Omega h^\kappa_{tx} |h^\kappa_x|_\delta\phi dxd\tau \non\\
   &\to& - 2\int_0^t\int_\Omega h^\kappa_{tx} |h^\kappa_x| \phi dxd\tau\ {\rm as}\ \delta\to 0 \non\\
   &=& - \int_0^t\int_\Omega \left(|h^\kappa_x|h^\kappa_{x} \right)_{t} \phi dxd\tau \non\\
   &=& -   \left(\left(|h^\kappa_x|h^\kappa_{x} \right)_{t}, \phi \right)_{Q_t}.
  \label{est-B5}
\end{eqnarray}

It then follows from \eq{est-B5} and \eq{est-B4} that
\begin{eqnarray}
 \left|\left(\left(|h^\kappa_x|h^\kappa_{x} \right)_{t}, \phi \right)_{Q_t} \right|
 = |I|
 \le  C\| \phi\|_{L^\infty(0,t;H^2_{\rm per}(\Omega))}.
  \label{est-B6}
\end{eqnarray}
Since $L^1(0,t;H^{-2}_{\rm per}(\Omega))$ is isometrically imbedded into  the dual space of $L^\infty(0,t;H^2_{\rm per}(\Omega))$, we complete the proof of Lemma~\ref{lm-est-B}.

\vskip0.2cm
We are going to study the asymptotic behavior of solution $h_B$ as
$B$ goes to zero. For this purpose we also need the following lemma.
\begin{Lemma}\label{lm5.2} Let $(0, T_e)\times\Omega$
be an open set in $\cR^+\times \cR^n$. Suppose functions $g_n, g$
are in $L^q((0, T_e)\times\Omega)$ for any given
$1 < q < \infty$, which satisfy
$$
 \|g_n\|_{L^q((0, T_e)\times\Omega)}\le C,\
 g_n \to g\ {\ a.e.\ in}\ (0, T_e)\times\Omega.
$$
Then $g_n$ converges to $g$ weakly in $L^q((0, T_e)\times\Omega)$.

\end{Lemma}

\vskip0.25cm
\noindent{\bf Proof of Theorem~\ref{thm1.1B}.}
Applying now the generalized case ($p_1=1$), of Aubin-Lions Lemma, i.e., Lemma~\ref{thm2.3} to the series $|(h_B)_x|(h_B)_x$, we assert from Lemma~\ref{lm-est-B} and Corollary~\ref{cor2.1} that
$$
 |(h_B)_x|(h_B)_x  \in L^\frac43(0,T_e;W^{1,\frac43}_{\rm per}(\Omega));\quad  (| (h_B)_x|(h_B)_x)_t\in L^1(0,T_e;H^{-2}_{\rm per}(\Omega)).
$$
This suggests us to choose
$$
 p_0=\frac43,\ p_1=1;\quad   B_0 = W^{1,\frac43}_{\rm per}(\Omega),
 \ B=L^2(\Omega),\ B_1 =  H^{-2}_{\rm per}(\Omega).
$$
We thus have $B_0\subset\subset B$,
$$
 |(h_B)_x|(h_B)_x \in L^{p_0}(0,T_e;B_0);\quad (| (h_B)_x|(h_B)_x)_t\in L^1(0,T_e;B_1),
$$
and conclude that
$
 |(h_B)_x|(h_B)_x
$ is compact in $L^\frac43(0,T_e;B)$.
  Hence we can select a subsequence, denote
  it by $h_{B_n}$, of $h_B$,  such that
$$
 |(h_{B_n})_x|(h_{B_n})_x\to \chi,\ a.e.,\ (t,x)\in\ Q_{T_e},
$$
$B_n\to 0$, as  $ n\to \infty$. It is easy to
see that the function $F: y \mapsto |y|y$
  is reversible, we obtain $(h_{B_n})_x\to F^{-1}(\chi) $  as  $ n\to 0$. By uniqueness of weak limit, we assert that $F^{-1}(\chi) = h_x $.

Recalling that $h_B$ satisfies
\begin{eqnarray}
(h_B,\varphi_{t})_{Q_{T_{e}}} &-& \frac{\alpha_{1}}{2} \Big( |(h_B)_{x}|(h_B)_{x} + 2B (h_B)_{x},\varphi_{x} \Big)_{Q_{T_{e}}} \nonumber\\
 &=& \Big((\alpha_{2}\sigma_i + \alpha_{3} )(|(h_B)_{x}| + B)  , \varphi\Big) _{Q_{T_{e}}}  - (h_{0},\varphi(0))_{\Omega} ,
 \label{1.11B1}
\end{eqnarray}
we need only study the limits of the most difficult terms, i.e.,
 the nonlinear terms like $(|(h_{B})_x|(h_{B})_x,\phi)_{Q_t}$.

Employing Lemma~\ref{lm5.2} we can easily pass the nonlinear terms
 to their limits. Thus
$h$ is a solution, in the sense of Definition~\ref{def5.1},
to   problem \eq{1.1} -- \eq{1.3} with $B=0$. And the proof of Theorem~\ref{thm1.1B} is thus complete.

\appendix
%

\section{Singular integrals}
For the sake of the readers' convenience, we include the
following theorem on the boundedness of singular integrals.
\begin{Theorem}[p. 48, Ref.~\cite{18}]\label{thm2.2}
Let $n\in\cN$ and $x\in\cR^n$. Suppose that the kernel $K(x)$ satisfies the conditions
\begin{eqnarray}
 |K(x)|\le C|x|^{-n},\quad  {\rm\ for\ all}\ |x|>0,
 \label{0.1}\\
 \int_{\{|x|\ge 2|y|\}}|K(x-y)-K(x)|dx\le C,\quad {\rm\ for\ all}\ |y|>0,
 \label{0.2}
\end{eqnarray}
and
\begin{eqnarray}
 \int_{\{R_{1}<|x|<R_{2}\} }K(x)dx=0,\quad 0<R_{1}<R_{2}<\infty,
 \label{0.3}
\end{eqnarray}
where $C$ is a positive constant. Let $ 1<p<\infty $,   for $f\in L^{p}(\cR^{n})$ we define
\begin{eqnarray}
 T_{\varepsilon}(f)(x)=\int_{\{|y|\ge \varepsilon\} }f(x-y)K(y)dy,\quad \varepsilon>0.
 \label{0.4}
\end{eqnarray}
Then there holds
\begin{eqnarray}
 \|T_{\varepsilon}(f)\|_{p}\le  C_{p}\|f\|_{p}
 \label{0.5}
\end{eqnarray}
here $C_{p}$ is a constant that is independent of $f$ and $\varepsilon$.
Also for each $f\in L^{p}(\cR^{n}),  \lim\limits_{\varepsilon\to 0}T_{\varepsilon}(f)=T(f)$
exists in $L^{p}$ norm. The operator $T$ so defined also satisfies the inequality \eqref{0.5}.
\end{Theorem}

The cancellation property alluded to is contained in condition \eqref{0.3}. This
hypothesis, together with \eqref{0.1}, \eqref{0.2}, allows us to prove the
$L^{2}$ boundedness and from this the $L^{p}$ convergence of the truncated
integrals \eqref{0.4}.

\vskip0.5cm
\noindent {\bf Acknowledgements.}
  The first author of this article is supported in part by Science and
  Technology Commission of Shanghai Municipality (Grant No. 20JC1413600); The
  corresponding author is supported in part  by the Hong Kong Research
  Grants Council General Research Fund 16302818 and
  Collaborative Research Fund C1005-19G.


\end{document}